\documentclass[11pt]{amsart}
\usepackage{amsmath,amssymb, graphicx, amscd,latexsym,here}
\makeatletter
\newtheorem{Theorem}{Theorem}
\newtheorem{Lemma}[Theorem]{Lemma}
\newtheorem{Corollary}[Theorem]{Corollary}
\newtheorem{Proposition}[Theorem]{Proposition}

\newtheorem{Remark}[Theorem]{Remark}

\newtheorem{Example}[Theorem]{Example}

\newtheorem{Conjecture}[Theorem]{Conjecture}

\newcommand{\eps}{\varepsilon}

\newcommand\vphi{\varphi}

\newcommand\al{\alpha}

\newcommand\be{\beta}

\newcommand\ga{\gamma}

\newcommand\cM{\mathcal  M}

\newcommand\cN{\mathcal  N}

\newcommand\BC{ {\mathbb C}}

\newcommand\BZ{{\mathbb  Z}}
\newcommand\BR{ {\mathbb  R}}
\newcommand\BP{ {\mathbb  P}}

\newcommand\bfu{\mbox {\bf  u}}

\newcommand\bfw{\mbox {\bf  w}}

\newcommand\bfz{\mbox {\bf  z}}

\newcommand\bfa{\mbox {\bf  a}}

\newcommand\nin{\noindent}

\newcommand\wtl{\widetilde}

\newcommand\id{\rm{id}}

\newcommand\codim{\rm{codim}\/}

\newcommand\lkn{{\rm{lkn}\/}}

\newcommand\rdeg{{\rm{rdeg}\/}}
\newcommand\pdeg{{\rm{pdeg}\/}}

\newcommand\inv{^{-1}}

\def\mapright#1{\smash{\mathop{\longrightarrow}\limits^{{#1}}}}

\def\mapdown#1{\Big\downarrow\rlap{$\vcenter{\hbox{$#1$}}$}}

\def\mapup#1{\Big\uparrow\rlap{$\vcenter{\hbox{$#1$}}$}}

\def\inv{^{-1}}

\begin{document}
\title[On mixed projective curves
]
{On mixed projective curves
}

\author
[M. Oka ]
{Mutsuo Oka }
\address{\vtop{
\hbox{Department of Mathematics}
\hbox{Tokyo  University of Science}
\hbox{26 Wakamiya-cho, Shinjuku-ku}
\hbox{Tokyo 162-8601}
\hbox{\rm{E-mail}: {\rm oka@rs.kagu.tus.ac.jp}}
}}

\begin{abstract}
Let $f(\bfz,\bar\bfz)$ be a mixed polar homogeneous polynomial
of $n$ variables $\bfz=(z_1,\dots, z_n)$. It
 defines a projective real algebraic  variety 
$V:=\{[\bfz]\in \BC\BP^{n-1}\,|\,f(\bfz,\bar\bfz)=0 \}$
in the projective space
$\BC\BP^{n-1}$. The behavior is  different from that of
 the projective hypersurface.
The topology is not uniquely determined by the degree of
 the variety even if $V$ is non-singular. We study a basic property of
 such a variety.
\end{abstract}
\maketitle

\maketitle

\section{Introduction}
Let $f(\bfz,\bar\bfz)$ be a polar weighted  homogeneous mixed polynomial 
with $\bfz=(z_1,\dots, z_n)\in \BC^n$.
Namely there exist  integers $(q_1,\dots, q_n)$ and 
$(p_1,\dots, p_n)$ and positive integers $d_r,\,d_p$ such that 
\begin{eqnarray*}\begin{split}
&f(t\circ \bfz,t\circ \bar\bfz)=t^{d_r}f(\bfz,\bar\bfz),\quad t\circ \bfz=(t^{q_1}z_1,\dots, t^{q_n}z_n),\,t\in \BR^+\\
&f(\rho\circ \bfz,\overline{\rho\circ\bfz})=\rho^{d_p}f(\bfz,\bar\bfz),\quad
\rho\circ \bfz=(\rho^{p_1}z_1,\dots, \rho^{p_n}z_n),\,\rho\in  \BC,\,|\rho|=1.
\end{split}\end{eqnarray*}
This gives $\BR^+\times S^1$ action by
\[
 (t,\rho)\circ \bfz=(t^{q_1}\rho^{p_1}z_1,\dots,
 t^{q_n}\rho^{p_n}z_n),\quad t\rho\in \BR^+\times S^1.
\]
The integers $d_r$ and $d_p$ are called the radial and the polar
degree respectively and we denote them as
$d_r=\rdeg\,f$ and $d_p=\pdeg\,f$.

We say that $f(\bfz,\bar \bfz)$ is {\em  strongly polar weighted  homogeneous}
if  $p_j=q_j$ for $j=1,\dots, n$. Then the associated $\BR^+\times S^1$
action on $\BC^n$ is in fact the $\BC^*$ action which is defined by
\[
 (\bfz,\tau)=((z_1,\dots, n),\tau)\mapsto \tau\circ\bfz=(z_1
 \tau^{p_1},\dots, z_n\tau^{p_n}),\,\,\tau\in \BC^*.
\]

We say $f(\bfz,\bar \bfz)$ is {\em  strongly polar  homogeneous}
if  further the weights satisfies the equalities
 $q_j=p_j=1$ for any $j$.
%
A strongly polar weighted homogeneous polynomial
$f(\bfz,\bar\bfz)$
satisfies the equality:
\begin{eqnarray}\label{action}
 f((t,\rho)\circ\bfz,\overline{(t,\rho)\circ\bfz})=t^{d_r}\rho^{d_p}f(\bfz,\bar\bfz),\quad
  (t,\rho)\in \BR^+\times S^1.
\end{eqnarray}

Assume that $f(\bfz,\bar \bfz)$ is a strongly polar weighted homogeneous polynomial
of radial degree $d_r$ and of polar degree $d_p$ respectively
and let $P=(p_1,\dots, p_n)$ be the weight vector.
Let $\tilde V$ be the mixed affine hypersurface
\[
 \tilde V=f\inv(0)=\{\bfz\in \BC^n\,|\, f(\bfz,\bar\bfz)=0\}.
\]
Let $\varphi:S^{2n-1}\setminus K \to S^1$  be the Milnor fibration
with  $K=\tilde V\cap S^{2n-1}$ and let $F$ be the fiber.
 Recall that
$\vphi(\bfz)=f(\bfz,\bar\bfz)/|f(\bfz,\bar\bfz)|$.
 Thus $F$ is defined by
\[
 F=\vphi\inv(1)=\{\bfz\in S^{2n-1}\setminus K\,|\, f(\bfz,\bar\bfz)>0\}
\]
We can equivalently consider the global fibration
$f:\BC^n-\tilde V\to \BC^*$. Then the Milnor fiber is identified with
the hypersurface
$f\inv(1)$.
The monodromy map $h:F\to F$ (in either case) is defined by
\[
 h(\bfz)=(\exp(\frac{2p_1\pi\,i}{d_p})z_1,\dots,\exp(\frac{2p_n\pi\,i}{d_p})z_n). 
\]
We consider also the wighted projective hypersurface $V$ defined by
\[
 V=\{(z_1:z_2:\cdots:z_n)\in \BC\BP(P)^{n-1}\,|\, f(\bfz,\bar\bfz)=0\}
\]
where $\BC\BP(P)^{n-1}$ is the weighted projective space defined by
the equivalence induced by the above $\BC^*$
action:
\[
 \bfz\sim \bfw\iff \exists \tau\in \BC^*,\, \bfw=\tau\circ \bfz.
\]
It is well-known that $\BC\BP(P)^{n-1}$ is an orbifold with at most  cyclic
quotient singularities.

By (\ref{action}),
 $\bfz\in f\inv(0)$  and $\bfz'\sim \bfz$, then $\bfz'\in f\inv(0)$.
Thus the hypersurface $V=\{[\bfz]\in \BC\BP^{n-1}(P)|f(\bfz)=0\}$ is well-defined.
Consider the quotient map
$\pi: S^{2n-1}\to \BC\BP(P)^{n-1}$ or $\pi: \BC^n\setminus\{O\}\to \BC\BP(P)^{n-1}$.
For the brevity's sake, we
 denote  the restrictions $\pi|F:F\to \BC\BP^{n-1}\setminus V$
and $\pi|K: K\to V$ by the same $\pi$.
We are interested in the topology of $V$ and  the relation with the
Milnor fibration. 

In this paper,  we consider only the case of strictly polar homogeneous
polynomials.
It turns out that the topology of the smooth projective mixed
hypersurface $V$ is not an invariant of the degree $d_r,d_p$. 
However we will show that 
{\em the degree of $V$ is equal to the polar degree $d_p$}
(Theorem \ref{degree}).

In \S 5, we study the case $n=3$.
In this case, let $g$ be the genus of the mixed curve $V$ and put $q=d_p$. Then
it is known that the following inequality holds
(known as the Thom's conjecture and proved by Kronheimer-Mrowka, \cite{Kronheimer-Mrowka}).
\[
 g\ge \frac{(q-1)(q-2)}2.
\]
 We also 
 give examples of mixed projective curves in
$\BC\BP^2$ which shows that $g$ is unbounded when $q$ is fixed.

This is a continuation of our papers \cite{OkaMix,OkaBrieskorn} and we use the same 
notations.
\section{Milnor fibration and the Hopf fibration}

\subsection{Canonical orientation}
It is well known that a complex analytic smooth variety has a canonical 
orientation which comes from the complex structure (see for example p.18,
\cite{GriffithsEtc}).
Let $\tilde V=f\inv(0)$ be a mixed hypersurface.
Take a point $\bfa\in \tilde V$. We say that $\bfa$ is a {\em mixed singular}
 point of
$\tilde V$,
 if $\bfa$ is a critical point of the
mapping
$f:\BC^n\to \BC$.
Otherwise, $\bfa$ is a {\em mixed regular point}.
Note that a point $\bfa\in \tilde V$ to be a regular point as a point of a real analytic variety
is a necessary condition but not  a sufficient condition
for the regularity as a point on a mixed variety.
Recall that $\bfa$ is a mixed singular point if and only if
$df_{\bfa}:T_{\bfa}\BC^n\to T_{f(\bfa)}\BC$ is surjective.
This is equivalent to the existence of a
complex number $\al\in  \BC$ with $|\al|=1$ such that 
\[\begin{split}
 &\overline{df}(\bfa,\bar\bfa)=\al \bar d\,f(\bfa,\bar\bfa)\quad
\text{i.e.},\quad
\frac{\overline{\partial f}}{\partial z_j}(\bfa,\bar\bfa)=
\al \frac{{\partial f}}{\partial \bar z_j}(\bfa,\bar\bfa),\,j=1,\dots,
   n
\end{split}
\](\cite{OkaPolar}).
We assert that 
\begin{Proposition} There is a canonical orientation on the smooth part of a mixed hypersurface.
\end{Proposition}
\begin{proof} Take a regular point $\bfa\in \tilde V$.
The normal bundle $\cN$ of $\tilde V\subset \BC^n$ has a canonical orientation
 so that 
$df_{\bfa}: \cN_{\bfa}\to T_{f(\bfa,\bar\bfa)}\BC$ is an orientation
 preserving isomorphism.
This gives a canonical orientation on $\tilde V$ so that 
the ordered union of the oriented frames
$\{v_1,\dots, v_{2n-2},n_1,n_2\}$ of $T_{\bfa} \BC^n$ is the orientation of
 $\BC^n$ if and only if 
$\{v_1,\dots, v_{2n-2}\}$ is an oriented frame of 
$T_{\bfa}\tilde V$ where $\{n_1,n_2\}$ is an oriented frame of normal vectors.
\end{proof}
Consider a mixed homogeneous hypersurface
$\tilde V$ and let $V$ be the corresponding mixed projective
hypersurface.
\begin{Proposition}
Let $\bfa\in \tilde V\setminus \{O\}$. Then $\bfa\in \tilde V$ is a mixed
singular point of $\tilde V$ if and only if 
$\pi(\bfa)\in V$ is a mixed singular point.
\end{Proposition}
\begin{proof}
Assume that $\bfa=(a_1,\dots, a_n)\in \wtl V$ and $a_n\ne 0$ for simplicity.
Let $u_j=z_j/z_n,\,1\le j\le n-1$ be the affine coordinates of the chart
$U_n=\{z_n\ne 0\}$. Then $V\cap U_n$ is defined by
\[
 V\cap U_n=\{\bfu\in \BC^{n-1}\,|\, g(\bfu,\bar\bfu)=0\}
\]
where $\bfu=(u_1,\dots, u_{n-1})$ and 
$g(\bfu,\bar\bfu)$ is defined by
\[
 g(\bfu,\bar\bfu):=f({\bfu'},\overline{{\bfu'}}),\quad
 {\bfu'}:=(\bfu,1).
\]
Putting $q+2r=\rdeg\,f$ and $q=\pdeg\,f$, we observe that 
\begin{eqnarray}\label{section}
g(\bfu,\bar\bfu)=f(\bfz,\bar\bfz)/(z_n^{q+r}\bar z_n^r).
\end{eqnarray} Here $d_r-d_p=2r$ and $d_p=q$.
Write $a_n=r_n\exp(\theta_n i)$ in the polar coordinate.
Consider the hyperplane section $\wtl U_n=\BC^n\cap\{z_n=a_n\}$ and $\wtl f:=f|_{\wtl U_n}$.
Then we have the commutative diagram:
\[
 \begin{matrix}
\wtl U_n&\mapright{\wtl f}& \BC\\
\mapdown{\pi}&&\mapup{\beta}\\
U_n=\BC^{n-1}&\mapright{g}&\BC
\end{matrix}
\]
where $\beta$ is the multiplication with 
$r_n^{d_r}\exp(id_p\theta_n)$. 
This follows from the equality (\ref{section}).
Put $\al:=\pi(a)\in U_n\cap V$.
Then the above diagram says that 
$d\wtl f_a:T_a\wtl U_n\to T_{O}\BC$ is surjective if and only if
$dg_{\al}:T_\al U_n\to T_{O}\BC$ is surjective.
On the other hand,
$T_a\BC^n$ is a direct sum of $T_a\wtl U_n$ and the tangent space of the
 $\BR^+\times S^1$ orbit at $a$ and the latter space is in the kernel of
$df_a:T_a\BC^n\to T_{O}\BC$, as $\wtl V$ is invariant by the
 $\BR^+\times S^1$-action.
This shows that the surjectivities of two tangential maps
$df_a:T_a\BC^n\to T_{O}\BC$ and $dg_\al U_n\to T_{O}\BC$
are equivalent. Thus $a\in \wtl V$ is mixed non-singular if and only if
$\al\in V$ is mixed nonsingular.
\end{proof}
Now we consider the canonical orientation of $V$.
First we recall that the orientation of $\BC^*$ is given
by the frame $\{\frac{\partial}{\partial r}, \frac{\partial}{\partial
\theta}\}$
where $(r,\theta)$ is the polar coordinates of $\BC^*$.
The orientation of $\BC\BP^{n-1}$  as a complex manifold and 
the orientation of $\BC\BP^{n-1}$ coming from the Hopf bundle
using the local bundle structure $U_j\times \BC^*$ is the same.
Using the orientation of the affine hypersurface $\widetilde  V$
and the local product structure of the restriction of  the Hopf bundle
over $V$, we have a canonical orientation on (the smooth part of ) $V$.
One can easily see that the orientation as 
the local mixed hypersurface
$V\cap U_n=g\inv(0)\subset U_n$ is the same with the above orientation.

\subsection{Milnor fiber} Consider the Hopf fibration
$\pi: S^{2n-1}\to \BC\BP^{n-1}$ and its restriction to the Milnor fiber
$F=\{\bfz\in S^{2n-1}\,|\,f(\bfz,\bar\bfz)>0\}$. As $f$ is polar
weighted,
it is easy to see that 
$\pi:F\to \BC\BP^{n-1}\setminus V$ is a cyclic covering of order
$d_p$ and the 
covering transformation is generated by the monodromy map
\[
 h:\,F\to F,\quad \bfz\mapsto \exp(\frac{2\pi i}{d_p})\cdot \bfz.
\]
Thus we have 
\begin{Proposition}
 \begin{enumerate}\label{formula1}
\item
$\chi(F)=d_p\chi(\BC\BP^{n-1}\setminus V)$.
     \item  $\chi(\BC\BP^{n-1}\setminus V)=n-\chi(V)$ and 
     $\chi(V)=n-\chi(F)/d_p$.
\item We have the following exact sequence.
$$1\to \pi_1(F)\mapright{\pi_{\sharp}} \pi_1(\BC\BP^{n-1}\setminus V)\to
     \BZ/d_p\BZ\to 1.$$
 \end{enumerate}
\end{Proposition}
\begin{Corollary}
If $d_p=1$, the projection $\pi:F\to \BC\BP^{n-1}\setminus V$ is a
 diffeomorphism.
\end{Corollary}
The monodromy map $h:F\to F$ gives free $\BZ/d_p\BZ$ action on $F$. Thus
using the periodic monodromy argument in \cite{Milnor}, we get
\begin{Proposition}
The zeta function of $h:F\to F$ is given by
\[
 \zeta(t)=(1-t^{d_p})^{-\chi(F)/d_p}.
\]In particular, if $d_p=1$,
$h=\id_F$ and 
$\zeta(t)=(1-t)^{-\chi(F)}$.
\end{Proposition}
\section{Topology of mixed projective hypersurface}
We are interested in the topology of the mixed projective hypersurface.
Assume that $f(\bfz,\bar\bfz)$ is a mixed homogeneous polynomial of
radial degree $d_r$ and of polar degree $q$. 
Let $V$  be the corresponding projective hypersurface
$V=\{f(\bfz,\bar\bfz)=0\}\subset \BC\BP^{n-1}$.
In the case of smooth  complex algebraic hypersurfaces, the topology of 
$F$ or $V$ are determined by the degree $q$. In the case of mixed 
hypersurfaces, we will see later that 
the degree $q$ do not determine the topology of
the Milnor fibering of $f$ or the topology of $V$.
\subsection{Isolated singularity case}
We consider a mixed strongly polar homogeneous polynomial
 $f(\bfz,\bar\bfz)$ of polar degree $q$ and we assume  that $\wtl V=f\inv (0)$
has an isolated mixed  singularity at the origin.

Assume that {\em $F=f\inv(1)$ 
 has a homotopy type of a
bouquet
of spheres of dimension $n-1$}.
\begin{Proposition}\label{homology}
 Under the above assumption, the projective mixed hypersurface has the
 following  homology groups:
 \[
    H_j(V)=\begin{cases}
            \BZ,\quad & j:even, \, j\le 2(n-2),\,j\ne n-2\\
            0,\quad &j:odd,\, j<2(n-2),\,j\ne n-2.
            \end{cases}
 \]
 The middle homology group $H_{n-2}(V)$ is determined by 
$\chi(F)$ up to the torsion part.
 \end{Proposition}
 \begin{proof}First, by the assumption,
$\tilde H_j(F)=0$ for $j\le n-2$ and 
  $H_{n-1}(F)=\BZ^{\mu}$ where
  \[
   \mu=(-1)^{n-1}(\chi(F)-1).
  \]
 By the Wang sequence of the Milnor fibration
\[
\dots \to H_{j+1}(S^{2n-1}\setminus K)\to H_j(F)\mapright{h_*-\id}H_j(F)\to
H_j(S^{2n-1}\setminus K)\to\dots
\] we see that
  \[
   H_j(S^{2n-1}\setminus K)=0,\quad j<n-2
  \]
As for $H_{n-2}(S^{2n-1}\setminus K)$,
it is isomorphic to the cokernel of
$h_*-\id:H_{n-1}(F)\to H_{n-1}(F)$ and the rank of this cokernel
can be computed from the zeta function and the characteristic
  polynomial
$P_{n-1}(t)$ which are related by
\[
 (t^{d_p}-1)^{\chi(F)/d_p}=\zeta(t)=\frac{P_{n-1}(t)^{(-1)^{n-2}}}{(t-1)}.
\]
We leave the calculation to the reader.
  By the Alexander duality, we get
  \[
  \tilde H_j(K)=0,\quad  j<n-2.
  \]
  Now the assertion follows from the Gysin sequence of the Hopf
  fibration
  $\pi: K\to V$:
  \[
   \dots \to H_j(K)\to H_j(V)\to H_{j-2}(V)\to H_{j-1}(K)\to \dots
  \]
The argument is exactly same as that  for a projective hypersurface (see \cite{MR52:6016})
  \end{proof}
\begin{Remark} Assume that a mixed function $f(\bfz,\bar\bfz)$ is strongly 
non-degenerate. It is an open problem if
(i) $F$ is $(n-2)$-connected, and (ii) $F$ has a homotopy type of
 CW-complex of dimension $n-1$.
\end{Remark}
\subsection{Solutions and points in $\BC\BP^1$} Let us consider the case
$n=2$.
Let  $\cM(q+2r,q;2)$ be the set of the mixed polar homogeneous
polynomials
with 
the radial degree $q+2r$ and the polar degree $q$.
Let $f(\bfz,\bar\bfz)$ be a non-degenerate strongly polar homogeneous
polynomial 
in $\cM(q+2r,q;2)$ 
where $\bfz=(z_1,z_2)$. 
For brevity, we assume that $q,r> 0$.
 We are interested to compute the number of points
$V:=\{f(\bfz,\bar\bfz)=0\}$. This number is equal to the number of
link components of $S^3\cap \wtl V$ where 
$\wtl V=f\inv(0)\subset \BC^2$ and we denoted this number by $\lkn
(\wtl V)$ in \cite{OkaMix}.
In general, $f$ takes the following form:
\[\begin{split}
 f(\bfz,\bar\bfz)&=\sum_{\nu,\mu}c_{\nu,\mu}z_1^{\nu_1}z_2^{\nu_2}\bar
 z_1^{\mu_1}\bar z_2^{\mu_2}\end{split}
\]
where 
$ \nu_1+\nu_2=\mu_1+\mu_2+q $   and $ \mu_1+\mu_2=r$.
We may assume that there are no points of $V$ with $z_2=0$.
Thus we may consider the coordinate chart $\{z_2\ne 0\}$ with
$z=z_1/z_2$
as the coordinate.
To know the exact number of points of $V$, we need to know the number
 of complex solutions of the mixed polynomial:
\[\begin{split}
& \{z\in \BC\,|\, \sum_{\nu_1,\mu_1}c_{\nu_1,\mu_1}'z^{\nu_1}\bar
 z^{\mu_1}=0\},
\quad \nu_1+\mu_1\le q+2r,\, \mu_1\le r\\
\end{split}
\]
where $c_{\nu_1,\mu_1}'=\sum_{\nu_2,\mu_2}c_{\nu,\mu}$.
In fact, the number of solutions is not so easy to be computed as in the
case of complex polynomials.
\begin{Example}{\rm
Consider the equation:
\[
-2 z^2\bar z+t z^2+1=0,\,\,t\in \BC.
\]
This example is considered in  Example 59 of our previous paper \cite{OkaMix}.
We can see that 
for a ``small'' $t$, we have only one solution.
For example $t=0$, $z=\frac{1}{\root{3}\of 2}$.
For a ``large'' $t$, we have three solutions.
(For real numbers, $t$ is ``small'' if $-3<t<1$.) For example, put $t=3$.
Then we get 
$z=a$ and $1/9\pm \sqrt{26}i/9$ where
$a$ is the real root of $-2a^3+3a^2+1=0$.
}
\end{Example}
This example tells us that the number of solutions depends on the coefficients.
 However we have the following observation.
\begin{Proposition}
Assume that $f(\bfz,\bar\bfz)\in \cM(q+2r,q;2)$ and 
 let  $V=\{[\bfz]\in \BC\BP^1\,|\,f(\bfz,\bar\bfz)=0\}$
and let $F=f\inv(1)\subset \BC^2$.
Then $\al:=\sharp V$ can take (at least)
$q,q+2,\dots, q+2r$. The corresponding Euler characteristic of $F$
are $\chi(F)=q(2-\al)$.
\end{Proposition}
\begin{proof}
We consider the basic two strongly polar homogeneous polynomials:
\[\begin{split}
& f_{q,j}:=z_1^{q+j}\bar z_1^{j}+z_2^{q+j}\bar z_2^j\in \cM(q+2j,q;2)\\
& k_{\ell}:=(z_1^\ell-\be z_2^\ell)(\bar z_1^\ell-\ga \bar z_2^\ell)\in
  \cM(2\ell,0;2),\,\,\be,\,\ga\in \BC^*\end{split}
\]
By Theorem 10, \cite{OkaPolar}, $f_{q,j}$ is strongly polar homogeneous
and $\lkn(\wtl V(f_{q,j}))=q$. $k_\ell$ is obviously strongly polar
 homogeneous
of degree 0 and $\lkn(\wtl(k_\ell))=2\ell$.
Thus $f_{q,j}k_{r-j}$ for $0\le j\le r$ is strongly polar homogeneous polynomial
in $\cM(q+2r,q;2)$
and non-degenerate as long as $V(f_{q,j})\cap V(k_\ell)=\emptyset$.
As $2\ell+ (q+2j)=q+2r$,
 $\lkn(\wtl V(f_{q,j}k_{r-j}))=q+2(r-j)$ with $j=1,\dots, r$.
The latter assertion follows from Lemma 64 of \cite{OkaMix}.
\end{proof}
\begin{Conjecture} The possible $\al$ for 
$f\in \cM(q+2r,q;2)$ is exactly $\{q,q+2,\dots, q+2r\}$.
\end{Conjecture}
\section{Degree of mixed projective hypersurfaces}
Suppose that $f(\bfz,\bar\bfz)\in \cM(q+2r,q;n)$ be a strongly polar homogeneous
polynomial
and let
\[
 V=\{\bfz\in \BC\BP^{n-1}\,|\,f(\bfz,\bar\bfz)=0\}.
\]
 We assume that 
the singular locus $\Sigma V$ of
$V$ is either empty or 
$\codim_{\BR}\Sigma V\ge 2$.
We have observed that $V\setminus\Sigma V\subset \BC\BP^{n-1}$
is canonically oriented so that the union of the oriented frame of $T_PV$,
say
$\{v_1,\dots, v_{2(n-2)}\}$ and the frame of normal bundle
$\{w_1,w_2\}$ which is compatible with the local defining
 complex function $g_j$ on the affine chart $U_j=\{z_j\ne 0\}$ is the
oriented frame of
$\BC\BP^{n-1}$. 
(Recall that $g_j$ is a mixed function of the 
variables
$u_i=z_i/z_j,\,i\ne j$ defined by
$g_j(\bfu,\bar\bfu)=f(\bfz,\bar\bfz)/z_j^{q+r}/\bar z_j^r$.)
Thus it has a fundamental class $[V]\in H_{2n-4}(V;\BZ)$ by
Borel-Haefliger \cite{Borel-Haefliger}. The topological degree  of
$V$
is the integer $d$ so that $\iota_{*}[V]=d[\BC\BP^{n-2}]$
where $\iota:V\to \BC\BP^{n-1}$ is the inclusion map
and $[\BC\BP^{n-2}]$ is the homology class of a canonical hyperplane
$\BC\BP^{n-2}$.

\vspace{.3cm}
The main result of this paper is: 
\begin{Theorem}\label{degree}
The topological degree of $V$ is
equal to the polar degree $q$. Namely
the fundamental class $[V]$ corresponds to $q[\BC\BP^{n-2}]\in
 H_{2(n-2)}(\BC\BP^{n-1})$ by the inclusion mapping $\iota_*$.
\end{Theorem}
\begin{proof}
Suppose that $f$ is a non-degenerate mixed polynomial in
 $\cM(q+2r,q;n)$.
Take a generic 1-dimensional complex line $L$ which is 
 isomorphic to $\BC\BP^1$. Then the degree is given by the intersection
 number 
$[V]\cdot [L]$.
Now, changing the coordinates if necessary, we may assume that 
\begin{eqnarray}\label{eq1}
 L:\,z_j=a_{j1}z_1+a_{j2}z_2,\,\,j=3,\dots,n.
\end{eqnarray}
Substituting (\ref{eq1}) in $f(\bfz,\bar\bfz)$ to eliminate the
 variables
$z_3,\dots, z_n$, we see that 
the intersection $V\cap L$ are described by
\[
 g(z_1,z_2,\bar z_1,\bar z_2)=0,\quad
[z_1:z_2]\in L=\BC\BP^1.
\]
As $g$ is still polar homogeneous in $z_1,z_2$ under the restriction to $L$,
 $g$ is written as 
\[
  g(\bfw,\bar \bfw)=f(\bfw,\bar\bfw)|_L=\sum_{\nu,\mu}c_{\nu,\mu}z_1^{\mu_1}z_2^{\nu_2}
\bar z_1^{\mu_1}\bar z_2^{\mu_2},\,
\bfw=(z_1,z_2),\,
\]
where the summation are for the multi-integers
$\nu=(\nu_1,\nu_2),\mu=(\mu_1,\mu_2)$ such that 
\[\begin{split}
 |\nu|+&|\mu|=q+2r,\,|\nu|-|\mu|=q,\,
|\nu|=\nu_1+\nu_2,\,|\mu|=\mu_1+\mu_2.
\end{split}
\]
Thus the polynomial $g(z_1,z_2,\bar z_1,\bar z_2)$
 is a polar homogeneous polynomial of polar degree $q$.
Taking a linear change of coordinates if necessary,
we may assume that the intersections are in the affine space $z_2\ne 0$.
This implies that $g$ has a monomial $z_1^{q+r}\bar z_1^r$ with 
a non-zero coefficient.
Use the affine coordinate $w=z_1/z_2$ for the affine coordinate
 $\{z_2\ne 0\}$. Then $g$ takes the form:
\[
 g(w,\bar w)=c_0w^{q+r}+c_1w^{q+r-1}+\cdots+c_{q+r}
\]
where $c_j$ is a polynomial in $\bar w$ such that 
$ \deg_{\bar w}c_j\le r$
and by the assumption, we have  that $c_0=\sum_{i=0}^r c_{0i}\bar w^i$
with $c_{0r}\ne 0$.
Let $\{\al_1,\dots, \al_m\}=\{w\,|\,g(\al,\bar\al)=0\}$.
We can see easily that 
\[
 I(V,L;\al_j)=\frac 1{2\pi}\int_{|w-\al_j|=\eps}{\rm Gauss}(g)d\theta
\] where $w-\al_j=\eps\exp(i\theta)$ and 
 ${\rm{Gauss}}(g)(w,\bar w)=\theta'$ with $\theta'=\arg\,(g(w,\bar w))$
and $\eps$ is  a sufficiently small  positive number.
In fact, the orientation of $V$ is defined so that 
a frame $\{v_1,\dots, v_{2n-4}\}$ at $\al_j$ is positive if and only if 
$\{v_1,\dots, v_{2n-4},n_1,n_2\}$ are positive where $n_1,n_2$ are
 frames of the normal bundle of $V$ oriented by $f$.
On the other hand, $\{\frac{\partial}{\partial
 x},\frac{\partial}{\partial y}\}$
is also a frame of the normal bundle where $w=x+iy$.
The orientation 
$\{n_1,n_2\}$ and $\{\frac{\partial}{\partial
 x},\frac{\partial}{\partial y}\}$ are compatible if and only if
the Gauss map at $\al_j$ has the positive rotation.

Topologically the intersection number is the mapping degree of
  the Gauss mapping, considered as:
\[
 {\rm Gauss}(g): \{|w-\al_j|=\eps\}\cong S^1\to S^1.
\]
Take  a sufficiently large positive
 number $R$. Then by a standard argument, we see that
\[
 \sum_{j=1}^m\frac 1{2\pi}\int_{|w-\al_j|=\eps}
{\rm Gauss}(g)d\theta
=
\frac 1{2\pi}\int_{|w|=R}{\rm Gauss}(g)d\theta.
\]
The right hand  side is equal to
the mapping degree of
\[
 {\rm Gauss(g)}:\{|w|=R\}\cong S^1\to S^1
\] 
which is equal to $q$ by the next Lemma which completes the proof.
\end{proof}

\subsection{Residue formula for a monic  mixed polynomial}
Let $g(w,\bar w)=\sum_{a,b}c_{a,b}w^a\bar w^b$ be a mixed polynomial.
Put $d=\max\{a+b\,|\, c_{a,b}\ne 0\}$ and we call $d$ the radial degree
of $g$. We say that $g$ is {\em a monic mixed polynomial of degree $d$} if
$g$ has a unique monomial of radial degree $d$. 
\begin{Lemma}Assume that $g(w,\bar w)$ is a monic mixed polynomial of
 degree $d$ which is written as 
\begin{eqnarray*}  
\begin{split}
g(w,\bar w)&=c_0(\bar w)w^{q+r}+c_1(\bar w)w^{q+r-1}+\cdots+c_{q+r}\\
&c_j(\bar w)\in \BC[\bar w],\, \rdeg_{\bar w}c_j\le r,\, j=0,\dots, q+r\\
&c_0(\bar w)=c_{0r}\bar w^r+\cdots+c_{00},\, c_{0r}\ne 0
\end{split}
\end{eqnarray*}
 with $d=q+2r$.
Then 
\[
 \frac 1{2\pi}\int_{|w|=R}{\rm Gauss}(g)d\theta=q.
\]
\end{Lemma}
\begin{proof}
Consider the family
$g_t(w,\bar w)=(1-t)\,g(w,\bar w)+t\,h(w,\bar w)$ with
 $h(w,\bar w):=c_{0r}\,w^{q+r}\bar w^r$.
For a sufficiently large $R$, this gives a homotopy of the two Gauss maps
of $g|_{|w|=R}$ and $h(w,\bar w)|_{|w|=R}$.
The rotation number of the Gauss map
$h(w,\bar w)|_{|w|=R}$ is obviously $q$.
This proves the assertion.
\end{proof}

\section{Mixed projective curves}
In this section, we study basic examples in the projective surface
$\BC\BP^2$.
Thus  we assume that $n=3$.
\subsection{Milnor fibers}
Let $f(\bfz,\bar\bfz)$ be a strongly polar
weighted  homogeneous polynomial in three variables $\bfz=(z_1,z_2,z_3)$.
Let $F=f\inv(1)\subset \BC^3$ be the Milnor
fiber.
\begin{Proposition}\label{Euler}
Assume that $f$ is 1-convenient (see \cite{OkaMix} for the definition)
 non-degenerate, polar weighted homogeneous polynomial with isolated mixed singularities at the
 origin
and we assume that $f$ is  either  of a join type or
of a simplicial type which are described  below.
\[
 \begin{split}
&\text{Join type}:\,\,f_1(\bfz,\bar\bfz)=g(z_1,z_2,\bar z_1,\bar z_2)+z_3^{a_3+b_3}\bar
 z_3^{b_3},\\
&\text{Simplicial type:}\begin{cases}
&f_2(\bfz,\bar\bfz)=z_1^{a_1+b_1}\bar z_1^{b_1}z_2+z_2^{a_2+b_2}\bar z_2^{b_2}z_3
+z_3^{a_3+b_3}\bar z_3^{b_3}\\
&f_2'(\bfz,\bar\bfz)=z_1^{a_1+b_1}\bar z_1^{b_1}\bar
 z_2+z_2^{a_2+b_2}\bar z_2^{b_2}\bar z_3
+z_3^{a_3+b_3}\bar z_3^{b_3}\\
&f_3(\bfz,\bar\bfz)=z_1^{a_1+b_1}\bar z_1^{b_1}z_2+z_2^{a_2+b_2}\bar z_2^{b_2}z_3
+z_3^{a_3+b_3}\bar z_3^{b_3}z_1\\
&f_3'(\bfz,\bar\bfz)=z_1^{a_1+b_1}\bar z_1^{b_1}\bar
 z_2+z_2^{a_2+b_2}\bar z_2^{b_2}\bar z_3
+z_3^{a_3+b_3}\bar z_3^{b_3}\bar z_1
\end{cases}
\end{split}
\]
with $a_i,b_i>0,\,i=1,2,3$
where $g(z_1,z_2,\bar z_1,\bar z_2)$ is a  non-degenerate 
polar weighted homogeneous polynomial.
 Then the Milnor fibers $F(f_i),\,i=1,\dots, 3$ and $F(f_i'),\,i=2,3$ 
are  simply connected and 
they have  homotopy types of  bouquets of spheres $S^2\vee\dots \vee
 S^2$.
Let $\mu_g$ be the Milnor number of $g$.
The Euler characteristics and the Milnor numbers  are given as follows.
\[
 \begin{split}
&\chi(F(f_1))=(a_3-1)\mu_g+ 1,\quad \mu(f_1)=(a_3-1)\mu_g\\
&\chi(F(f_2))=\chi(F(f_2'))=a_1a_2a_3-a_2a_3+a_3,\\
 &\qquad\qquad\qquad \mu(f_2)=\mu(f_2')=\chi(F(f_2))-1,\\
&\chi(F(f_3))=a_1a_2a_3+1,\quad \mu(f_3)=a_1a_2a_3,\\
&\chi(F(f_3'))=a_1a_2a_3-1\quad \mu(f_3')=a_1a_2a_3-2,\\
\end{split}
\]
\end{Proposition}
\begin{proof} 
We consider  first  $F_1=f_1\inv(1)$
where
\[
 f_1(\bfz,\bar\bfz)=g(z_1,z_2,\bar z_1,\bar z_2)+z_3^{a_3+b_3}\bar z_3^{b_3}
\]
where $g(z_1,z_2,\bar z_1,\bar z_2)$ is a convenient  non-degenerate
polar weighted homogeneous
 polynomial.
For two variables case, the Milnor fiber $F_g$
of $g(z_1,z_2)$ has the homotopy type of a bouquet of $S^1$ as
it is a connected open Riemann surface (Proposition 38, \cite{OkaMix}). Let $\mu_g$ be the Milnor number (=the first
 Betti number) of $F_g$. Then the Milnor fiber $F_1$ of $f_1$  is
 homotopic to the join $F_g*\Omega_a$
where $\Omega_a$ is the set of $a$-th roots of unity (\cite{Molina}).
This join is obviously homotopic to a bouquet of $\mu_g(a-1)$ $S^2$
 spheres.

Consider $F_2=f_2\inv(1)$
or $F_2'=f_2\inv(1)$.
The Euler characteristic  can be easily
 computed from the additivity of the Euler characteristics,
applied on  the toric
stratification
\[
 F_2=F_2^{*\{1,2,3\}}\amalg F_2^{*\{2,3\}}\amalg F_2^{*\{3\}},\quad
 F_2'=F_2'^{*\{1,2,3\}}\amalg F_2'^{*\{2,3\}}\amalg F_2'^{*\{3\}},\quad
\]
and Theorem 10 of \cite{OkaPolar},
where $F_2^{*I}$ is defined by $F_2\cap \BC^{*I}$ and 
 \[
  \BC^{*I}=\{\bfz\in \BC^3|z_i\ne 0,\,i\in I,\, z_j=0,\,j\notin I\}\]
for $I\subset \{1,2,3\}$.

The Euler characteristics of $F_3=f_3\inv(1)$ and $F_3'={f_3'}\inv(1)$
can be computed in the exact same way.
The assertion on the homotopy types are now obtained simultaneously as
 follows.
First $F_j^{*\{1,2,3\}}$ and ${F_j'}^{*\{1,2,3\}}$ 
are CW-complex of the  dimension 2 by Theorem 10, \cite{OkaPolar}.
The secondly $F_j$ and  ${F_j'}$  are simply connected by the
 1-convenience
assumption (\cite{OkaPolar}). Using the above decomposition and
 Mayer-Vietoris exact sequences,
we see that the (reduced) homology groups are non-trivial only 
on the dimension 
 2 and no torsion on $H_2(F_j)$
and $H_2({F_j'})$ for $j=2,3$.
Thus by the  Whitehead theorem (see for example \cite{Spanier}), 
we conclude that $F_j$  and ${F_j'}$  are homotopic to
 bouquets of two dimensional spheres.
\end{proof}
\subsection{Projective mixed curves.}
We consider projective curves of degree $q$:
\[
 C=\{[z_1:z_2:z_3]\in\BC\BP^2\,|\, f(z_1,z_2,z_3)=0\}
\]
where $f$ is a strongly polar homogeneous polynomial with
 $\pdeg\,f=q$.
We have seen that the topological degree of $C$ is $q$
by Theorem \ref{degree}. The genus $g$ of $C$ is not
 an invariant of $q$. 
Recall that for a differentiable curve $C$ of genus $g$, embedded in
$\BC\BP^2$,
with the
topological degree $q$, we have  the following Thom's inequality, which
was conjectured by Thom and proved by for example
Kronheimer-Mrowka \cite{Kronheimer-Mrowka}:
\[
 g\ge \frac{(q-1)(q-2)}2
\]
where the right side number is the genus of algebraic curves of degree
$q$, given by 
the Pl\"ucker formula.
Recall that for a mixed strongly polar homogeneous polynomial,
the genus and the Euler characteristic
of the Milnor fiber are related as follows ( by (2) of Proposition 3): 
\begin{Proposition}\label{genus}
\begin{eqnarray}\label{formula1}
 2-2g=3-\frac{\chi(F)}{q}
\end{eqnarray}
where 
\[
 F=\{(z_1,z_2,z_3)\in\BC^3\,|\, f(z_1,z_2,z_3,\bar z_1,\bar z_2,\bar z_3)=1\}.
\]
\end{Proposition}
Now we will see some examples which shows that $\chi(F)$ is not an
invariant of $q$.

\nin
{\bf I. Simplicial polynomials.}
We consider the following simplicial polar homogeneous polynomials
of polar degree $q$.
\[
 \begin{split}
&f_{s_1}(\bfz,\bar\bfz)=z_1^{q+r}\bar z_1^r+z_2^{q+r}\bar
z_2^r+z_3^{q+r}\bar z_3^r\\
&f_{s_2}(\bfz,\bar\bfz)=z_1^{q+r-1}\bar z_1^rz_2+z_2^{q+r-1}\bar
z_2^rz_3+z_3^{q+r}\bar z_3^r\\
&f_{s_3}(\bfz,\bar\bfz)=z_1^{q+r-1}\bar z_1^rz_2+z_2^{q+r-1}\bar
z_2^rz_3+z_3^{q+r-1}\bar z_3^rz_1\\
&f_{s_4}(\bfz,\bar\bfz)=z_1^{q+r+1}\bar z_1^r\bar z_2+z_2^{q+r+1}\bar
z_2^r\bar z_3+z_3^{q+r}\bar z_3^r\\
&f_{s_5}(\bfz,\bar\bfz)=z_1^{q+r+1}\bar z_1^r\bar z_2+z_2^{q+r+1}\bar
z_2^r\bar z_3+z_3^{q+r+1}\bar z_3^r\bar z_1\\
\end{split}
\]
Let $F_{s_i}$ be the Milnor fiber of $f_{s_i}$ 
and let $C_{s_i}$ be the corresponding projective curves for $i=1,\dots,
5$.
First, the Euler characteristic of the Milnor fibers
and the genera are given by Proposition \ref{Euler}
and Proposition \ref{genus} as follows.
\[
 \begin{split}
&\chi(F_{s_i})=q^3-3q^2+3q,\quad
  g(C_{s_i})=\frac{(q-1)(q-2)}2,\,i=1,2,3\\
&\chi(F_{s_4})=q(q^2+q+1),\quad g(C_{s_4})=\frac {q(q+1)}2\\
&\chi(F_{s_5})=q(q^2+3q+3),\quad g(C_{s_4})=\frac {(q+2)(q+1)}2\\
\end{split}
\]
In \cite{OkaBrieskorn},
we have shown that  $C_{s_1}$ and $C_{s_2}$ are isomorphic to  algebraic
plane curves defined by
the associated homogeneous polynomials of degree $q$:
\begin{eqnarray*}
 &g_{s_1}(\bfz)=z_1^{q-1}z_2+z_2^{q-1}z_3+z_3^{q}\\
&g_{s_2}(\bfz)=z_1^{q-1}z_2+z_2^{q-1}z_3+z_3^{q}.
\end{eqnarray*}
We also expect that $C_{s_3}$ is isotopic to the algebraic 
curve
\[
 z_1^{q-1}z_2+z_2^{q-1}z_3+z_3^{q-1}z_1=0,
\]
as the genus of $C_{s_3}$ suggests it (see also \cite{OkaBrieskorn}).

\vspace{.2cm}
\nin
{\bf II.} We consider the following join type polar homogeneous
polynomial.
\[\begin{split}
&h_{j}(\bfz,\bar\bfz)= g_j(\bfw,\bar\bfw)+z_3^{q+r}\bar z_3^r,\\
   &g_j(\bfw,\bar \bfw)=(w_1^{q+j}\bar w_1^j+w_2^{q+j}\bar
   w_2^{j})(w_1^{r-j}-\al w_2^{r-j})(\bar w_1^{r-j}-\be \bar
   w_2^{r-j}),\,\\
 &\qquad \qquad  0\le j\le r.
   \end{split}
\] ($\al,\be\in \BC^*$ are generic.)
The the Milnor fiber $F_{g_j}$ of $g_j$ is connected.
The Euler characteristic of $\chi(F_{g_j}^*)$ ($F_{g_j}^*=F_{g_j}\cap \BC^{*2}$)
is  given by  $\chi(F_{g_j}^*)=-r_g\,q$ where $r_g$ is the link component number
of $g=0$ which is  $q+2(r-j)$. Thus
 $\chi(F_{g_j})=\chi(F_{g_j}^*)+2q$ where the last terms come from
$\chi(F_{g_j})^{*I})$ with $I=\{1\}$ or $\{2\}$. Thus
$\chi(F_{g_j})=q(q-2+2(r-j))$ and
$g=\frac{(q-1)(q-2+2(r-j))}2$.
We observe that 
 the genus can take the following values by taking $j=r,\dots, 0$:
\[
 \frac{(q-1)(q-2)}2,\,\frac{(q-1)q}2,\,\dots, \frac{(q-1)(q+2r-2)}2.
\]
As we can take the positive number $r$ arbitrary large,
we have that
\begin{Proposition} There exist differentiable curves
embedded in $\BC\BP^2$ with a fixed degree
$q\ge 2$ whose genera are given as:
\[
 \{g_0+k(q-1)\,|\,k=0,1,\dots\},\quad g_0:=\frac{(q-1)(q-2)}2.
\]
\end{Proposition}

In particular, taking $q=2$, we obtain:
\begin{Corollary}
 For any smooth  surface $S$ of genus $g$, there is an embedding $S\subset
 \BC\BP^2$
 so that the degree of $S$ is 2.
 \end{Corollary}
\bibliographystyle{abbrv}
\def\cprime{$'$} \def\cprime{$'$} \def\cprime{$'$} \def\cprime{$'$}
  \def\cprime{$'$} \def\cprime{$'$} \def\cprime{$'$} \def\cprime{$'$}

\end{document}